# DOUBLING AND PROJECTION: A METHOD OF CONSTRUCTING TWO-LEVEL DESIGNS OF RESOLUTION IV[1]

By Hegang H. Chen and Ching-Shui Cheng

*University of Maryland School of Medicine, and Academia Sinica and University of California, Berkeley*

Given a two-level regular fractional factorial design of resolution IV, the method of doubling produces another design of resolution IV which doubles both the run size and the number of factors of the initial design. On the other hand, the projection of a design of resolution IV onto a subset of factors is of resolution IV or higher. Recent work in the literature of projective geometry essentially determines the structures of all regular designs of resolution IV with $n \geq N/4 + 1$ in terms of doubling and projection, where $N$ is the run size and $n$ is the number of factors. These results imply that, for instance, all regular designs of resolution IV with $5N/16 < n \leq N/2$ must be projections of the regular design of resolution IV with $N/2$ factors. We show that, for $9N/32 \leq n \leq 5N/16$, all minimum aberration designs are projections of the design with $5N/16$ factors which is constructed by repeatedly doubling the $2^{5-1}$ design defined by $I = ABCDE$. To prove this result, we also derive some properties of doubling, including an identity that relates the wordlength pattern of a design to that of its double and a result that does the same for the alias patterns of two-factor interactions.

**1. Introduction.** Doubling is a simple but powerful method of constructing two-level fractional factorial designs, in particular, those of resolution IV. Suppose $\mathbf{X}$ is an $N \times n$ matrix with two distinct entries, 1 and $-1$. Then the *double* of $\mathbf{X}$, denoted by $D(\mathbf{X})$, is the $2N \times 2n$ matrix

$$\begin{bmatrix} \mathbf{X} & \mathbf{X} \\ \mathbf{X} & -\mathbf{X} \end{bmatrix},$$

Received June 2004; revised March 2005.
[1]Supported in part by NSF Grant DMS-00-71438.

that is,

$$D(\mathbf{X}) = \begin{bmatrix} 1 & 1 \\ 1 & -1 \end{bmatrix} \otimes \mathbf{X},$$

where $\otimes$ is the Kronecker product. Suppose $\mathbf{X}$ defines an $N$-run design for $n$ two-level factors, where the two levels are denoted by 1 and $-1$, each column of $\mathbf{X}$ corresponds to a factor and each row of $\mathbf{X}$ defines a factor–level combination. Then $D(\mathbf{X})$ defines a design which doubles both the run size and the number of factors of $\mathbf{X}$. Such a method was used by Plackett and Burman [10] in their classical paper on orthogonal main-effect plans.

An $N \times m$ submatrix of $\mathbf{X}$, where $m \leq n$, is called a *projection* of $\mathbf{X}$ (onto $m$ factors). Equivalently, such a design can be obtained by deleting $n - m$ columns (factors) from $\mathbf{X}$. We allow the possibility $m = n$, so that a design is considered as its own projection onto all columns. If $\mathbf{X}$ is of resolution IV, then all projections of $\mathbf{X}$ are of resolution IV or higher. A less obvious fact is that if $\mathbf{X}$ is of resolution IV, then $D(\mathbf{X})$ is also of resolution IV. Thus, starting with a design of resolution IV, the combined operation of doubling followed by projection yields a design of resolution IV or higher.

Some recent results in the literature of finite projective geometry [2, 3, 8] essentially characterize, in terms of doubling and projection, the structures of *regular* designs of resolution IV with $n \geq N/4 + 1$. These elegant and difficult results have important implications in statistical design, but are not easily accessible to statisticians since they were stated in the language of projective geometry. One purpose of this paper is to review some of these results and rephrase them in design language. This should be useful to people interested in the statistical design of experiments. In the past few years we have encountered several occasions where some special cases of these results were either re-discovered by statisticians or would have been very helpful if they had been become more widely known in the statistical literature.

We also present some new results, including an identity that relates the wordlength pattern of a design to that of its double and a result that does the same for the alias patterns of two-factor interactions. These results and some other basic properties of doubling are presented in Section 2. Characterization of regular designs of resolution IV with $n \geq N/4 + 1$ in terms of doubling and projection is discussed in Section 3. Section 4 shows how the results in Sections 2 and 3 can be used to investigate certain minimum aberration designs. In particular, we show that, for $9N/32 \leq n \leq 5N/16$, a minimum aberration design can be obtained by deleting factors from the minimum aberration design with $5N/16$ factors which can be constructed by repeatedly doubling the $2^{5-1}$ design defined by $I = ABCDE$.

We conclude this section with a review of some basic terminology and notation. A factorial design is called an orthogonal array of strength $t$ if, in all projections onto $t$ factors, all factor–level combinations appear the same number of times. *Regular* fractional factorial designs are those which can be constructed by using defining relations and are discussed in many textbooks on experimental design; see, for example, [11]. Each interaction that appears in the defining relation is called a defining word, and the resolution of a regular design is defined as the length of the shortest defining word. It is well known that a regular design of resolution $R$ is an orthogonal array of strength $R-1$. For each positive integer $i$, let $B_i$ be the number of defining words of length $i$. Then the resolution is equal to the smallest $i$ such that $B_i > 0$. The sequence $(B_1, B_2, \ldots, B_n)$ is called the wordlength pattern of the design. The *minimum aberration* criterion introduced by Fries and Hunter [9] chooses a design by sequentially minimizing $B_1$, $B_2$, $B_3$, ....

Throughout this paper $N$ and $n$ denote the run size and the number of factors, respectively. We shall restrict consideration to two-level designs only. Under a regular two-level design, $N$ must be a power of 2, say, $N = 2^{n-p}$. This design is a $\frac{1}{2^p}$-fraction of a complete $2^n$ factorial, and is usually referred to as a $2^{n-p}$ design. It is well known that a regular two-level design of resolution III must have $n \leq N-1$, and a regular two-level design of resolution IV must have $n \leq N/2$. A design of resolution III is called saturated if $n = N-1$, and a design of resolution IV is called saturated if $n = N/2$.

**2. Some basic properties of doubling.** We first note that if $\mathbf{X}$ is a Hadamard matrix of order $N$, then $D(\mathbf{X})$ is a Hadamard matrix of order $2N$. In fact, this was what Plackett and Burman used in their 1946 classic paper. The following properties can easily be established:

THEOREM 2.1. *If $\mathbf{X}$ is an orthogonal array of strength two, then $D(\mathbf{X})$ is also an orthogonal array of strength two. Likewise, if $\mathbf{X}$ is an orthogonal array of strength three, then $D(\mathbf{X})$ is an orthogonal array of strength three.*

Theorem 2.1 has no counterpart for designs of higher strength. In fact, for two columns $\mathbf{a}$ and $\mathbf{b}$ of $\mathbf{X}$, $D(\mathbf{X})$ must have four columns of the form $\begin{bmatrix} \mathbf{a} & \mathbf{b} & \mathbf{a} & \mathbf{b} \\ \mathbf{a} & \mathbf{b} & -\mathbf{a} & -\mathbf{b} \end{bmatrix}$, whose componentwise product has all the entries equal to 1. Therefore, $D(\mathbf{X})$ cannot have strength higher than three. For regular designs, these four columns would lead to a defining word of length four.

Another important elementary fact is that saturated regular designs of resolution III are unique (up to isomorphism). Such a design of size $2^k$ can be obtained by deleting the first column of

$$\underbrace{\begin{bmatrix} 1 & 1 \\ 1 & -1 \end{bmatrix} \otimes \begin{bmatrix} 1 & 1 \\ 1 & -1 \end{bmatrix} \otimes \cdots \otimes \begin{bmatrix} 1 & 1 \\ 1 & -1 \end{bmatrix}}_{k}.$$

In other words, saturated regular designs of resolution III can be obtained by deleting a column of 1's after successively doubling the $2 \times 2$ matrix $\begin{bmatrix} 1 & 1 \\ 1 & -1 \end{bmatrix}$.

Since all regular designs of resolution III or higher can be constructed by deleting a subset of columns from a saturated regular design of resolution III (or, equivalently, by projecting the saturated regular design of resolution III onto the complementary set of columns), all two-level regular designs of resolution III or higher can be constructed by the operation of doubling followed by deletion (or projection). These observations also imply that if $\mathbf{X}$ is a regular design, then $D(\mathbf{X})$ is regular. In other words, in Theorem 2.1, if $\mathbf{X}$ is a regular design of resolution III (resp. IV), then $D(\mathbf{X})$ is also a regular design of resolution III (resp. IV). However, by the comment following Theorem 2.1, if $\mathbf{X}$ is a regular design of resolution higher than IV, then $D(\mathbf{X})$ is a regular design of resolution IV only.

Saturated regular designs of resolution IV are also unique (up to isomorphism). One way of constructing such designs is to *foldover* saturated regular designs of resolution III. We would like to present another neat and compact method of construction using doubling. Consider the $2^2$ complete factorial

$$\begin{bmatrix} 1 & 1 \\ 1 & -1 \\ -1 & 1 \\ -1 & -1 \end{bmatrix},$$

whose number of factors is precisely half of the run size. Successively doubling the $2^2$ complete factorial results in designs whose number of factors is half of the run size, and by the discussion at the end of the previous paragraph, are regular designs of resolution IV. This implies that saturated regular designs of resolution IV can be obtained by successively doubling the $2^2$ complete factorial.

As mentioned earlier, the designs obtained by doubling those of resolution III or higher are of resolution III or IV. To study statistical properties of such designs, it is important to consider the alias pattern of two-factor interactions. Suppose $N = 2^{n-p}$. Among the $2^n - 1$ factorial effects, $2^p - 1$ appear in the defining relation. The rest are divided into $g \equiv 2^{n-p} - 1$ alias sets, each of size $2^p$. Without loss of generality, assume that the first $f \equiv 2^{n-p} - 1 - n$ of these alias sets does not contain main effects. For each $1 \leq i \leq g$, let $m_i$ be the number of two-factor interactions in the $i$th alias set. Cheng, Steinberg and Sun [7] showed that

$$(2.1) \qquad B_3 = \tfrac{1}{3}\left[\binom{n}{2} - \sum_{i=1}^{f} m_i\right]$$

and

$$(2.2) \qquad B_4 = \tfrac{1}{6}\left[\sum_{i=1}^{g} m_i^2 - \binom{n}{2}\right].$$

Furthermore, estimation capacity, a measure of model robustness discussed in [7], can be expressed in terms of the $m_i$'s. The following result relates the $m_i$ values of a design to those of its double.

THEOREM 2.2. *Suppose a regular design* $\mathbf{X}$ *has* $g$ *alias sets, the first* $f$ *of which do not contain main effects, and for* $1 \leq i \leq g$, $m_i$ *is the number of two-factor interactions in the ith alias set. Let the corresponding numbers of* $D(\mathbf{X})$ *be* $g^*, f^*$ *and* $m_i^*$, *respectively. Then* $g^* = 2g+1, f^* = 2f+1$, $m_1^* = m_2^* = 2m_1, m_3^* = m_4^* = 2m_2, \ldots, m_{2f-1}^* = m_{2f}^* = 2m_f$, $m_{2f+1}^* = n$, *and* $m_{2f+2}^* = m_{2f+3}^* = 2m_{f+1}, \ldots, m_{g^*-1}^* = m_{g^*}^* = 2m_g.$

PROOF. Suppose $\mathbf{X}$ is a $2^{n-p}$ design. Then $D(\mathbf{X})$ is a $2^{2n-(n+p-1)}$ design, and $g = 2^{n-p} - 1$, $f = 2^{n-p} - 1 - n$, $g^* = 2^{2n-(n+p-1)} - 1$, $f^* = 2^{2n-(n+p-1)} - 2n - 1$. The relations $f^* = 2f+1$ and $g^* = 2g+1$ follow.

To each factor in $\mathbf{X}$, say $A$, there correspond two factors in $D(\mathbf{X})$. Suppose the column of $\mathbf{X}$ corresponding to $A$ is $\mathbf{a}$. We shall denote the factor in $D(\mathbf{X})$ corresponding to the column $\begin{bmatrix}\mathbf{a}\\\mathbf{a}\end{bmatrix}$ by $A^+$ and the factor corresponding to $\begin{bmatrix}\mathbf{a}\\-\mathbf{a}\end{bmatrix}$ by $A^-$. Then it is easy to see that, for any two factors $A$ and $B$ in $\mathbf{X}$, under $D(\mathbf{X})$ $A^+B^+$ and $A^-B^-$ are aliased, $A^+B^-$ and $A^-B^+$ are aliased, and $A^+A^-$ and $B^+B^-$ are aliased. Consequently, if the two-factor interactions $AB$ and $CD$ are aliased under $\mathbf{X}$, then $A^+B^+$, $A^-B^-$, $C^+D^+$ and $C^-D^-$ are aliased under $D(\mathbf{X})$, and $A^+B^-$, $A^-B^+$, $C^+D^-$ and $C^-D^+$ are also aliased. If the main effect $A$ and two-factor interaction $CD$ are aliased under $\mathbf{X}$, then $A^+$, $C^+D^+$ and $C^-D^-$ are aliased under $D(\mathbf{X})$, and $A^-$, $C^+D^-$ and $C^-D^+$ are also aliased. Using these facts, it can be seen that each alias set of two-factor interactions under $\mathbf{X}$ determines two alias sets of two-factor interactions under $D(\mathbf{X})$, each of which is twice as large as the original alias set under $\mathbf{X}$. Furthermore, all the $n$ two-factor interactions of the form $A^+A^-$ constitute another alias set which does not contain main effects. □

We have the following relationship between the wordlength pattern of $\mathbf{X}$ and that of $D(\mathbf{X})$.

THEOREM 2.3. *Let $B_k$ and $B_k^*$ be the number of defining words of length $k$ of $\mathbf{X}$ and $D(\mathbf{X})$, respectively. Then*

$$B_k^* = \begin{cases} \sum_{s=0}^{\min[(k-1)/2, n]} B_{k-2s} \cdot \binom{n-k+2s}{s} 2^{k-2s-1}, & \text{if } k \text{ is odd;} \\ \sum_{s=0}^{\min[k/2-1, n]} B_{k-2s} \cdot \binom{n-k+2s}{s} 2^{k-2s-1} + \binom{n}{\frac{k}{2}}, & \\ & \text{if } k \text{ is even and } k/2 \text{ is even;} \\ \sum_{s=0}^{\min[k/2-1, n]} B_{k-2s} \cdot \binom{n-k+2s}{s} 2^{k-2s-1}, & \text{if } k \text{ is even and } k/2 \text{ is odd.} \end{cases}$$

PROOF. If $A_{i_1} \cdots A_{i_k}$ is a defining word of $\mathbf{X}$, then $A_{i_1}^{j_1} \cdots A_{i_k}^{j_k}$, where each $j_l$ is a $+$ or $-$, is a defining word of $D(\mathbf{X})$ as long as an even number of the $j_l$'s are $-$'s. When $s$ is even, for any $s$ distinct factors $A_{i_1}, \ldots, A_{i_s}$, $A_{i_1}^+ A_{i_1}^- \cdots A_{i_s}^+ A_{i_s}^-$ is also a defining word of $D(\mathbf{X})$. In general, a defining word of $D(\mathbf{X})$ is of the form $A_{i_1}^+ A_{i_1}^- \cdots A_{i_s}^+ A_{i_s}^- A_{i_{s+1}}^{j_{s+1}} \cdots A_{i_{s+t}}^{j_{s+t}}$, where $A_{i_1}, \ldots, A_{i_{s+t}}$ are distinct factors of $\mathbf{X}$, $A_{i_{s+1}} \cdots A_{i_{s+t}}$ is a defining word of $\mathbf{X}$, and an even (resp. odd) number of the $j_l$'s are $-$'s when $s$ is even (resp. odd). If $A_{i_1}^+ A_{i_1}^- \cdots A_{i_s}^+ A_{i_s}^- A_{i_{s+1}}^{j_{s+1}} \cdots A_{i_{s+t}}^{j_{s+t}}$ is of length $k$, then $t = k - 2s$. This implies that $0 \leq s \leq \min[k/2, n]$. On the other hand, when $s$ is odd, we must have $t \geq 1$. Therefore, we have $s \leq \min[(k-1)/2, n]$ when $k$ is odd, $s \leq \min[k/2, n]$ when $k$ is even and $k/2$ is even, and $s \leq \min[k/2-1, n]$ when $k$ is even and $k/2$ is odd. The theorem then follows from the fact that, for given $s$, there are $2^{t-1} = 2^{k-2s-1}$ ways to choose $(j_{s+1}, \ldots, j_{s+t})$ if $t \geq 1$, $B_{k-2s}$ ways to choose $A_{i_{s+1}}, \ldots, A_{i_{s+t}}$ and $\binom{n-k+2s}{s}$ ways to choose $A_{i_1}, \ldots, A_{i_s}$. Note that the last term $\binom{n}{\frac{k}{2}}$ in the case where $k/2$ is even corresponds to $s = k/2$. In this case, there are $\binom{n}{\frac{k}{2}}$ ways to choose $A_{i_1}^+ A_{i_1}^- \cdots A_{i_{k/2}}^+ A_{i_{k/2}}^-$. □

In particular, for designs of resolution III or higher we have $B_3^* = 4B_3$ and $B_4^* = 8B_4 + \binom{n}{2}$. These identities can also be derived by using (2.1), (2.2) and Theorem 2.2.

The following is an immediate consequence of Theorem 2.3.

COROLLARY 2.4. *Given two regular designs $\mathbf{X}_1$ and $\mathbf{X}_2$, $\mathbf{X}_1$ has less aberration than $\mathbf{X}_2$ if and only if $D(\mathbf{X}_1)$ has less aberration than $D(\mathbf{X}_2)$.*

**3. Maximal designs of resolution IV.** One major difference between designs of resolution III and IV is that all regular designs of resolution III can be constructed by deleting factors from saturated regular designs of resolution III. On the other hand, not all regular designs of resolution IV can be obtained by deleting factors from saturated regular designs of resolution IV: while all designs obtained by deleting factors from saturated regular designs of resolution IV are the so-called even designs which have no defining words of odd lengths, there are designs of resolution IV which are not even designs.

We say that a regular design of resolution IV or higher is maximal if its resolution reduces to three whenever an extra factor is added. Clearly, the saturated regular designs of resolution IV are maximal. If a design is not maximal, then at least one factor can be added so that the design is still of resolution IV or higher. One can keep adding factors until it becomes maximal. Therefore, if a regular design of resolution IV is not maximal, then it can be obtained by deleting factors from a maximal design. Because of the importance of this fact, we state it formally as a proposition:

PROPOSITION 3.1. *Every regular design of resolution* IV *is a projection of a certain maximal regular design of resolution* IV *or higher.*

We can also define maximal regular designs of resolution III. It turns out that, for a given run size, there is only one maximal regular design of resolution III: the saturated regular design of resolution III. As mentioned earlier, saturated regular designs of resolution IV are maximal, but they are not the only maximal designs of resolution IV.

Theoretically speaking, for a fixed number of runs, if we can determine all the maximal designs of resolution IV or higher, then all regular designs of resolution IV can be constructed by projecting the maximal designs onto subsets of factors. Recently some significant progress in the determination of such maximal designs has been made in the literature of finite projective geometry. Note that maximal designs of resolution IV or higher are equivalent to maximal caps in a finite projective geometry.

We first state a simple but interesting characterization of maximal regular designs of resolution IV.

THEOREM 3.2. *A regular design of resolution* IV *is maximal if and only if* $m_i > 0$ *for all* $i = 1, \ldots, f$.

This result was first stated in the coding-theoretic language; see, for example, [2]. Chen and Cheng [5] rephrased it in the above form. They also defined the notion of *estimation index*. Another way to state the result in Theorem 3.2 is that a regular design of resolution IV is maximal if and only if its estimation index is equal to 2. Since we can estimate one effect from each alias set assuming that the other effects in the same alias set are negligible, Theorem 3.2 says that a regular design of resolution IV is maximal if and only if all the available degrees of freedom can be used to estimate main effects and two-factor interactions. Such designs are said to be *second-order saturated* by Block and Mee [1].

The following result, whose geometric version can also be found in [2], reveals the crucial role played by the method of doubling in constructing designs of resolution IV.

THEOREM 3.3. *Let* **X** *be a regular design of resolution* IV *or higher. Then* **X** *is maximal if and only if $D(\mathbf{X})$ is maximal.*

One can see that Theorem 3.3 follows immediately from Theorem 2.2 and Theorem 3.3.

The following two key results essentially determine the structures of regular resolution IV designs with $N/4 + 1 \leq n \leq N/2$.

THEOREM 3.4 ([2, 8]). *Every maximal regular design of resolution* IV *with $N/4+2 \leq n \leq N/2$ can be obtained by doubling a maximal regular design of resolution* IV *or higher.*

THEOREM 3.5 ([3]). *For each $N = 2^k$ with $k \geq 4$ there exists at least one maximal regular design of resolution* IV *or higher with $n = N/4 + 1$.*

One family of maximal designs with $n = N/4 + 1$ can be found in Tang, Ma, Ingram and Wang's [13] study of designs with maximum number of clear two-factor interactions.

A complete search shows that there are two maximal regular 16-run designs of resolution IV or higher. One is the saturated $2^{8-4}$ design of resolution IV, and the other is the $2^{5-1}$ design defined by $I = ABCDE$, whose existence is assured by Theorem 3.5. Repeatedly doubling the former yields all larger saturated regular designs of resolution IV, while successively doubling the latter (a design of resolution V) leads to a family of maximal designs of resolution IV with $n = 5N/16$, for $N = 2^k$, $k \geq 5$. Since there are no 16-run maximal regular designs of resolution IV with $5 < n < 8$, by Theorem 3.4, for all $N = 2^k$ with $k \geq 5$, there are no maximal designs of resolution IV with $5N/16 < n < N/2$. This leads to the following important conclusion.

COROLLARY 3.6. *For $5N/16 < n < N/2$, a regular design of resolution* IV *must be obtained by deleting columns from saturated designs of resolution* IV. *All such designs are even designs.*

Doubling the two maximal 16-run regular designs of resolution IV, we obtain the saturated $2^{16-11}$ design of resolution IV and a $2^{10-5}$ design, both of which are maximal. By Theorem 3.5 there exists at least one maximal regular 32-run design of resolution IV with $n = 9$. In fact, there is exactly one such design. Therefore, for $9 \leq n \leq 16$ there are exactly three 32-run maximal designs of resolution IV: a $2^{16-11}$, a $2^{10-5}$ and a $2^{9-4}$. The last one and its repeated doubles constitute a family of maximal regular designs of resolution IV with $n = 9N/32$, for $N = 2^k$, $k \geq 5$. Again, there are no maximal regular designs of resolution IV with $9N/32 < n < 5N/16$. Thus, for the $n$'s in this range, a regular design of resolution IV must be obtained by deleting columns from either the saturated design of resolution IV or the maximal regular design of resolution IV with $n = 5N/16$.

Doubling the three maximal 32-run designs of resolution IV, we obtain the saturated $2^{32-26}$ design of resolution IV, a $2^{20-14}$ design and a $2^{18-12}$ design, all of which are maximal. By Theorem 3.5 there exists at least one maximal regular 64-run design of resolution IV with $n = 17$. Block and Mee's [1] complete search shows that there are five such designs. Therefore, for $17 \leq n \leq 32$ there are eight 64-run maximal regular designs of resolution IV: a $2^{32-26}$, a $2^{20-14}$, a $2^{18-12}$ and five $2^{17-11}$. One can also conclude that there are exactly five maximal regular designs of resolution IV with $n = 17N/64$, for all $N = 2^k$ with $k \geq 6$.

Now it is clear that if $N = 2^k$, $k \geq 4$, then for $n \geq N/4 + 1$ a maximal regular design of resolution IV or higher must have

(3.1) $\qquad n \in \{N/2, 5N/16, 9N/32, 17N/64, 33N/128, \dots\}.$

Conversely, for each integer $n = (2^i + 1)N/2^{i+2}$, there exists at least one maximal regular $N$-run design of resolution IV or higher with $n$ factors. A maximal regular design of resolution IV or higher with $n = (2^i + 1)N/2^{i+2}$ and $N = 2^k$, $k > i + 2$, can be constructed by repeatedly doubling a maximal regular $2^{i+2}$-run design of resolution IV or higher with $2^i + 1$ factors.

**4. Some results on minimum aberration designs.** By the discussion in the previous section, there are maximal regular designs of resolution IV with $n = N/2$, $5N/16, 9N/32, 17N/64, 33N/128, \dots$. Those with $n = N/2$, the saturated regular designs of resolution IV, are known to have minimum aberration. For $5N/16 < n < N/2$, minimum aberration designs (in fact, all regular designs of resolution IV) must be projections of the saturated regular design of resolution IV. Butler [4] addressed the issue of deleting factors from the saturated regular design of resolution IV so that the resulting design has minimum aberration. This is reminiscent of the complementary design theory of Chen and Hedayat [6], Tang and Wu [14] and Suen, Chen and Wu [12] that deals with how to find a set of factors so that its complement in the saturated regular design of resolution III has minimum aberration.

In an unpublished work, N. A. Butler found the minimum aberration design with $n = 5N/16$ which, using the terminology in this paper, is the maximal regular design of resolution IV with $n = 5N/16$. In this section we shall show that, for $9N/32 \leq n \leq 5N/16$, the minimum aberration designs are projections of the maximal regular design of resolution IV with $n = 5N/16$. Thus, although for $n > 9N/32$ minimum aberration designs are projections of the maximal regular design of resolution IV with either $N/2$ or $5N/16$ factors, the first two in (3.1), the pattern breaks down at $9N/32$. Even the maximal design of resolution IV with $n = 9N/32$ itself does not have minimum aberration.

Before proceeding to the proof of this result, we shall explore a bit more the alias pattern of two-factor interactions under the maximal design of resolution IV with $n = 5N/16$. Suppose $N = 16 \cdot 2^t$, where $t \geq 1$, and let $\mathbf{X}^*$ be the maximal design of resolution IV with $5 \cdot 2^t$ factors. Then $\mathbf{X}^*$ can be obtained by doubling the $2^{5-1}$ design defined by $I = ABCDE$ $t$ times. Suppose $\mathbf{a}, \mathbf{b}, \mathbf{c}, \mathbf{d}$ and $\mathbf{e}$ are the five columns of this $2^{5-1}$ design corresponding to factors $A, B, C, D$ and $E$, respectively. Then each of $\mathbf{a}, \mathbf{b}, \mathbf{c}, \mathbf{d}$ and $\mathbf{e}$ generates $2^t$ columns of $\mathbf{X}^*$. For example, each of the $2^t$ columns of $\mathbf{X}^*$ generated by $\mathbf{a}$ takes the form $\mathbf{x}_t \otimes \cdots \otimes \mathbf{x}_1 \otimes \mathbf{a}$, where $\mathbf{x}_i = [1, 1]^T$ or $[1, -1]^T$. We shall denote the corresponding factor of $\mathbf{X}^*$ by $A^{\mathbf{j}}$, where $\mathbf{j} = (j_1, \ldots, j_t)$, with $j_i = 1$ if $\mathbf{x}_i = [1, 1]^T$ and $j_i = -1$ if $\mathbf{x}_i = [1, -1]^T$. Notation such as $B^{\mathbf{j}}, C^{\mathbf{j}}, D^{\mathbf{j}}$ and $E^{\mathbf{j}}$ is similarly defined. The $5 \cdot 2^t$ factors of $\mathbf{X}^*$ are thus partitioned into five groups each of size $2^t$.

Any two of the five factors $A, B, C, D$ and $E$, say $X$ and $Y$, generate $2^t \cdot 2^t$ two-factor interactions of the form $X^{\mathbf{i}} Y^{\mathbf{j}}$, where $\mathbf{i}$ and $\mathbf{j}$ are $1 \times t$ vectors with entries 1 or $-1$. These interactions, $\binom{5}{2} 2^t \cdot 2^t = 10 \cdot 2^t \cdot 2^t$ in total, are called between-group two-factor interactions. Those of the form $X^{\mathbf{i}} X^{\mathbf{j}}$ with $\mathbf{i} \neq \mathbf{j}$, $5\binom{2^t}{2} = 5 \cdot 2^{t-1}(2^t - 1)$ in total, are called *within-group* two-factor interactions.

Since the $2^{5-1}$ design defined by $I = ABCDE$ is of resolution V, there is exactly one two-factor interaction in each of its ten alias sets not containing main effects. By applying Theorem 2.2 $t$ times, we see that, under $\mathbf{X}^*$, there are $10 \cdot 2^t + (2^t - 1)$ alias sets not containing main effects, $10 \cdot 2^t$ of which each containing $2^t$ two-factor interactions, and $2^t - 1$ of which each containing $5 \cdot 2^{t-1}$ two-factor interactions.

It can be seen that the $2^t \cdot 2^t$ between-group two-factor interactions arising from the same pair $(X, Y)$ are distributed evenly in $2^t$ alias sets of size $2^t$. Each of these $2^t$ alias sets is of the form $\{X^{\mathbf{i}} Y^{\mathbf{j}} : \mathbf{i} \odot \mathbf{j} = \mathbf{k}\}$, where $\mathbf{k}$ is a $1 \times t$ vector with entries 1 or $-1$ and $\mathbf{i} \odot \mathbf{j}$ is the componentwise product of $\mathbf{i}$ and $\mathbf{j}$. We denote this alias set by $XY_{\mathbf{k}}$. The ten possible pairs of $X$ and $Y$ account for the $10 \cdot 2^t$ alias sets of size $2^t$ mentioned in the previous paragraph. On the other hand, the $2^{t-1}(2^t - 1)$ within-group two-factor interactions $X^{\mathbf{i}} X^{\mathbf{j}}$ arising from the same $X$ are distributed evenly in the remaining $2^t - 1$ alias sets. Each of these alias sets consists of the $5 \cdot 2^{t-1}$ interactions $A^{\mathbf{i}} A^{\mathbf{j}}$, $B^{\mathbf{i}} B^{\mathbf{j}}$, $C^{\mathbf{i}} C^{\mathbf{j}}$, $D^{\mathbf{i}} D^{\mathbf{j}}$, $E^{\mathbf{i}} E^{\mathbf{j}}$ with $\mathbf{i} \odot \mathbf{j} = \mathbf{k}$. (Note that $X^{\mathbf{i}} X^{\mathbf{j}}$ is the same as $X^{\mathbf{j}} X^{\mathbf{i}}$ and $\mathbf{k} \neq \mathbf{1}$, where $\mathbf{1}$ is the vector of 1's.) We denote this alias set by $W_{\mathbf{k}}$.

One key property that is important for the proof is that each of the $2^t$ factors of the form $X^{\mathbf{i}}$ appears in *exactly one* of the $2^t$ two-factor interactions in each alias set $XY_{\mathbf{k}}$. As a consequence, if $u$ factors generated by $X$ are deleted from $\mathbf{X}^*$, then the number of two-factor interactions in each of these $2^t$ alias sets is reduced by $u$. In this case, since each factor of the $2^{5-1}$ design can be coupled with four other factors, the number of two-factor interactions in $4 \cdot 2^t$ of the $10 \cdot 2^t$ alias sets of size $2^t$ is reduced to $2^t - u$; that in each of the other $6 \cdot 2^t$ alias sets remains to be $2^t$.

It can also be seen that each of the $5 \cdot 2^t$ factors appears in exactly one within-group two-factor interaction in each $W_{\mathbf{k}}$.

Now we are ready to prove the following theorem.

THEOREM 4.1. *For any $N = 2^k$, $k \geq 5$, and $9N/32 \leq n \leq 5N/16$, the minimum aberration design must be a projection of the design constructed by repeatedly doubling the $2^{5-1}$ design defined by $I = ABCDE$.*

PROOF. Let $N = 16 \cdot 2^t$ and $u = 5N/16 - n$. Then since $5N/16 - 9N/32 = 2^{t-1}$, we have $0 \leq u \leq 2^{t-1}$.

The minimum aberration design with $N = 32$ and $n = 9$ is known to have at least one zero among the $m_i$'s, $1 \leq i \leq f$ ([7], page 91), and therefore is not maximal. In fact, it is obtained by deleting one factor from the maximal design with ten factors. By repeatedly applying Corollary 2.4, we see that, for each $N = 2^k$ with $k \geq 5$, the maximal regular design of resolution IV with $9N/32$ factors does not have minimum aberration. Now since there are only two maximal regular designs of resolution IV with more than $9N/32$ factors, one with $N/2$ factors and the other with $5N/16$ factors, it is enough to show that, for $N = 16 \cdot 2^t$ with $t \geq 2$, there is at least one $(5N/16 - u)$-factor projection of the maximal regular design of resolution IV with $5N/16$ factors that has less aberration than *all* projections of the saturated regular design of resolution IV. As before, let $\mathbf{X}^*$ be the maximal regular design of resolution IV with $5N/16$ factors.

Now consider a design $\mathbf{X}$ obtained by deleting from $\mathbf{X}^*$ $u$ factors that are generated by the same factor $A$ of the $2^{5-1}$ design defined by $I = ABCDE$. We shall show that $\mathbf{X}$ has fewer defining words of length four or, equivalently, by (2.2), a smaller value of $\sum_{i=1}^{f} m_i^2$ than all projections of the saturated regular design of resolution IV.

By the observations preceding the statement of the current theorem, under $\mathbf{X}$, $4 \cdot 2^t$ of the $10 \cdot 2^t$ alias sets of between-group two-factor interactions have $m_i$ equal to $2^t - u$, and the other $6 \cdot 2^t$ alias sets have $m_i$ equal to $2^t$. Now we provide an upper bound on the sum of squares of the $m_i$ values over the alias sets of within-group two-factor interactions.

Under $\mathbf{X}^*$, in each of the $2^t - 1$ alias sets of within-group two-factor interactions, $2^{t-1}$ of the $5 \cdot 2^{t-1}$ two-factor interactions involve factors generated by $A$. Thus, when only factors generated by $A$ are deleted, for each of these alias sets, the resulting $m_i$ satisfies $4 \cdot 2^{t-1} \leq m_i \leq 5 \cdot 2^{t-1}$. Consequently, an upper bound on the sum of squares of the $m_i$ values over these alias sets can be obtained by making as many of the $m_i$ values equal to $5 \cdot 2^{t-1}$ or $4 \cdot 2^{t-1}$ as possible.

When $u$ factors generated by $A$ are deleted, $\binom{u}{2} + u(2^t - u)$ within-group two-factor interactions are also deleted. Write

$$u(u+1)/2 = a \cdot 2^{t-1} + b,$$

where $a$ and $b$ are nonnegative integers such that $b < 2^{t-1}$. Then

$$(4.1) \qquad \binom{u}{2} + u(2^t - u) = (2u - a)2^{t-1} - b.$$

By the observation two paragraphs above, an upper bound on the sum of squares of the $m_i$ values over the $2^t - 1$ alias sets of within-group two-factor interactions can be obtained by assuming that $2u - a - 1$ of the $m_i$ values are equal to $4 \cdot 2^{t-1}$, one is equal to $4 \cdot 2^{t-1} + b$, and the remaining $2^t - 2u + a - 1$ values are equal to $5 \cdot 2^{t-1}$. Combining this with the $m_i$ values for the alias sets of between-group interactions obtained earlier, we conclude that if $u$ factors generated by $A$ are deleted from $\mathbf{X}^*$, then

$$(4.2) \quad \sum_{i=1}^{f} m_i^2 \leq 4 \cdot 2^t \cdot (2^t - u)^2 + 6 \cdot 2^t \cdot (2^t)^2 + (4 \cdot 2^{t-1} + b)^2$$
$$+ (2^t - 2u + a - 1)(5 \cdot 2^{t-1})^2 + (2u - a - 1)(4 \cdot 2^{t-1})^2.$$

On the other hand, a saturated regular design of resolution IV has all the $\binom{5 \cdot 2^t - u}{2}$ two-factor interactions in $N/2 - 1 = 2^{t+3} - 1$ alias sets. Thus, a design obtained by deleting factors from a saturated regular design of resolution IV can have at most $2^{t+3} - 1$ nonzero $m_i$'s. For such a design,

$$(4.3) \qquad \sum_{i=1}^{f} m_i^2 \geq \left[ \frac{\binom{5 \cdot 2^t - u}{2}}{2^{t+3} - 1} \right]^2 \cdot (2^{t+3} - 1).$$

It is sufficient to show that the right side of (4.2) is less than that of (4.3). This can be verified by tedious calculations using (4.1) and the assumption that $0 \leq b < 2^{t-1}$. The details are omitted. $\square$

Note that the case $u = 0$ provides an alternative proof of Butler's result on the optimality of maximal designs of resolution IV with $n = 5N/16$.

Theorem 4.1 leaves open the issue of which projection of the maximal design of resolution IV with $n = 5N/16$ has minimum aberration. A complementary design theory in the same spirit as that for saturated regular designs of resolution III and IV needs to be developed.

Suppose $u$ factors are deleted from $\mathbf{X}^*$. For $X = A, B, C, D, E$, let $n_X$ be the number of factors deleted from those generated by $X$. Among the $n_X n_Y$ pairs $(\mathbf{i}, \mathbf{j})$ where $X^\mathbf{i}$ and $Y^\mathbf{j}$ are deleted, let $n_{X,Y}^\mathbf{k}$ be the total number such that $\mathbf{i} \odot \mathbf{j} = \mathbf{k}$; this is the number of between-group two-factor interactions formed by the deleted $X^\mathbf{i}$ and $Y^\mathbf{j}$'s that belong to the alias set $XY_\mathbf{k}$. Similarly, among the $\binom{n_X}{2}$ ordered pairs $(\mathbf{i}, \mathbf{j})$ where $X^\mathbf{i}$ and $X^\mathbf{j}$ are deleted, let $n_{X,X}^\mathbf{k}$ be the total number such that $\mathbf{i} \odot \mathbf{j} = \mathbf{k}$; this is the number of within-group two-factor interactions formed by the deleted $X^\mathbf{i}$'s that belong to the alias set $W_\mathbf{k}$. Then the number of two-factor interactions in $XY_\mathbf{k}$ is reduced by $n_X + n_Y - n_{X,Y}^\mathbf{k}$, and the number of two-factor interactions in $W_\mathbf{k}$ is reduced by $n_A + n_B + n_C + n_D + n_E - n_{A,A}^\mathbf{k} - n_{B,B}^\mathbf{k} - n_{C,C}^\mathbf{k} - n_{D,D}^\mathbf{k} - n_{E,E}^\mathbf{k} = u - n_{A,A}^\mathbf{k} - n_{B,B}^\mathbf{k} - n_{C,C}^\mathbf{k} - n_{D,D}^\mathbf{k} - n_{E,E}^\mathbf{k}$. Thus, a minimum aberration projection of $\mathbf{X}^*$ must minimize

$$\sum_{X,Y} \sum_\mathbf{k} (2^t - n_X - n_Y + n_{X,Y}^\mathbf{k})^2$$
$$+ \sum_{\mathbf{k} \neq \mathbf{1}} (5 \times 2^{t-1} - u + n_{A,A}^\mathbf{k} + n_{B,B}^\mathbf{k} + n_{C,C}^\mathbf{k} + n_{D,D}^\mathbf{k} + n_{E,E}^\mathbf{k})^2,$$

where, in the first term, the first sum is over the ten possible $(X,Y)$'s and the second sum is over all the $1 \times t$ vectors of 1's and $-1$'s. While this can be solved without difficulty for small $u$'s, more general results need to be developed. We expect an optimal strategy to delete the factors one at a time alternately from the five groups of factors generated by $A, B, C, D, E$, while making the two-factor interactions formed by the deleted factors as uniformly distributed among the alias sets as possible.

**Acknowledgments.** We would like to thank the referees for their helpful comments.

DIVISION OF BIOSTATISTICS AND BIOINFORMATICS
DEPARTMENT OF EPIDEMIOLOGY
 AND PREVENTIVE MEDICINE
UNIVERSITY OF MARYLAND SCHOOL OF MEDICINE
113 HOWARD HALL
660 WEST REDWOOD STREET
BALTIMORE, MARYLAND 21201
USA
E-MAIL: hchen@epi.umaryland.edu

DEPARTMENT OF STATISTICS
UNIVERSITY OF CALIFORNIA
BERKELEY, CALIFORNIA 94720
USA
E-MAIL: cheng@stat.berkeley.edu